\documentclass{article}
\usepackage{amsfonts}
\usepackage{amsmath}

\newtheorem{theorem}{Theorem}

\newenvironment{proof}[1][Proof]{\noindent\textbf{#1.} }{\ \rule{0.5em}{0.5em}}
\input{tcilatex}
\begin{document}

\title{Dimension Free Growth Results For Vector-Valued Functions Of Several
Complex Variables}
\author{Faruk F. Abi-Khuzam \\
%EndAName
The American University Of Beirut}
\maketitle

\begin{abstract}
Let $f$ be an entire function of finite order $\rho \in (0,1)$. The maximum
modulus $M(r)$ of $f$ and the counting function of the zeros $N(r)$ are
connected by the following best possible growth inequality known as
Valiron's Theorem:%
\begin{equation*}
\limsup_{r\rightarrow \infty }\frac{N(r)}{\log M(r)}\geq \frac{\pi \rho }{%
\sin \pi \rho }.
\end{equation*}%
For functions subharmonic in $%
%TCIMACRO{\U{211d} }%
%BeginExpansion
\mathbb{R}
%EndExpansion
^{d}$, Hayman obtained a corresponding result with a best possible constant
involving the dimension $d.$ For the special case of an entire function on $%
%TCIMACRO{\U{2102} }%
%BeginExpansion
\mathbb{C}
%EndExpansion
^{d\text{ }},$ we obtain a corresponding dimension-free, best possible
inequality.
\end{abstract}

\bigskip

\section{\protect\bigskip Introduction}

\bigskip

Let $f$ be an non-constant entire function. Let $M(r;f)$ $%
=\sup_{|z|=r}|f(z)|,$ the maximum modulus, and $N(r,\frac{1}{f})$ the
smoothed counting function for the zeros of $f.$ The order of $f$ is defined
by

\begin{equation}
\rho =\lim \sup_{r\rightarrow \infty }\frac{\log \log M(r;f)}{\log r}.
\end{equation}%
The question of comparing the growth of the two real functions $\log M(r;f)$
and $N(r,\frac{1}{f})$ at infinity, in particular supplying a sharp lower
bound in terms of $\rho $ for $\lim \sup_{r\rightarrow \infty }\frac{N(r,%
\frac{1}{f})}{\log M(r;f)}$ , is a rather difficult question in the theory
that remains largely unsolved at present. For recent developments see \cite%
{AS} and \cite{SW}. In the special case where $0<\rho <1$, a solution
was\bigskip\ supplied in the following well-known result of Polya-Valiron 
\cite{polya}, \cite{val1},\cite{val2} 
\begin{equation*}
\lim \sup_{r\rightarrow \infty }\frac{N(r,\frac{1}{f})}{\log M(r;f)}\geq 
\frac{\sin \pi \rho }{\pi \rho }.
\end{equation*}%
Equality holds for entire functions of order $\rho \in (0,1)$ whose zeros
are regularly distributed over one ray in the complex plane.

Similar sharp results are also known if $\log M(r;f)$ is replaced by $%
m_{p}(r;\log |f|)$, the $L^{p}$ mean defined by, 
\begin{equation*}
m_{p}(r;\log |f|)=\left\{ \frac{1}{2\pi }\int_{-\pi }^{\pi }|\log
|f(re^{i\theta })||^{p}d\theta \right\} ^{1/p},
\end{equation*}%
provided the order $\rho \in (0,1).$If we put 
\begin{equation*}
\psi (\theta )=\frac{\pi \rho }{\sin \pi \rho }\cos \rho \theta ,-\pi \leq
\theta \leq \pi ,
\end{equation*}%
then, for each $1\leq p<\infty $, we have the sharp inequality%
\begin{equation*}
\lim \sup_{r\rightarrow \infty }\frac{N(r,\frac{1}{f})}{m_{p}(r;\log |f|)}%
\geq \frac{1}{\left\{ \frac{1}{2\pi }\int_{-\pi }^{\pi }|\psi (\theta
)|^{p}d\theta \right\} ^{1/p}}.
\end{equation*}%
It should be pointed out that, in the special case $p=2,$ this later result
holds true for meromorphic functions of any finite order by a well-known
result of Miles-Shea.\bigskip

The purpose of this paper is to obtain a dimension-free analogue of the the
above results for entire functions of several complex variables, and of
order less than $1.$

\section{Statement of results}

Let $f:%
%TCIMACRO{\U{2102} }%
%BeginExpansion
\mathbb{C}
%EndExpansion
^{d}\rightarrow 
%TCIMACRO{\U{2102} }%
%BeginExpansion
\mathbb{C}
%EndExpansion
$ be an entire function and assume for simplicity that $f(0)=1$. Denote by $%
S=S^{2d-1}$ the unit sphere in $%
%TCIMACRO{\U{2102} }%
%BeginExpansion
\mathbb{C}
%EndExpansion
^{d}$. For $\zeta \in S$, the slice function $f_{\zeta }:%
%TCIMACRO{\U{2102} }%
%BeginExpansion
\mathbb{C}
%EndExpansion
\rightarrow 
%TCIMACRO{\U{2102} }%
%BeginExpansion
\mathbb{C}
%EndExpansion
$, is defined by 
\begin{equation*}
f_{\zeta }(z)=f(z\zeta )\text{.}
\end{equation*}%
Clearly, $f_{\zeta }$ is, for each fixed $\zeta ,$an entire function with $%
f_{\zeta }(0)=1$, and the real functions $M(r,f_{\zeta }),$ $N(r,\frac{1}{%
f_{\zeta }})$ $,$ and $m_{p}(r;\log |f_{\zeta }|)$ are well-defined for $%
r\geq 0$. We now define the functions $M(r;f),$ $n_{\max }(r),$ and $%
(m_{p})_{\max }$ by 
\begin{equation*}
M(r;f)=\sup_{\zeta \in S}M(r,f_{\zeta }),n_{\max }(r;f))=\sup_{\zeta \in
S}n(r,\frac{1}{f_{\zeta }}),(m_{p})_{\max }(r;f)=\sup_{\zeta \in
S}m_{p}(r;\log |f_{\zeta }|),
\end{equation*}%
and introduce the two "smoothed" functions%
\begin{equation*}
N_{\max }(r;f)=\sup_{\zeta \in S}N(r,\frac{1}{f_{\zeta }}),N(r,f)=%
\int_{0}^{r}\frac{n_{\max }(t;f))}{t}dt.
\end{equation*}%
The assumption $f(0)=1$ implies that $n_{\max }$ vanishes identically in a
neighbourhood of $0$. Also $n_{\max }$ is non-decreasing and so the function 
$N(r;f)$ is well-defined. The functions $\log M(r,f)$ and $N_{\max }(r;f)$
are increasing convex functions of $\log r$, and it is easily seen that $%
M(r;f)$ is the maximum modulus of $f$ on the sphere in $%
%TCIMACRO{\U{2102} }%
%BeginExpansion
\mathbb{C}
%EndExpansion
^{d}$ centered at $O$ and having radius $r$. But $N(r;f)\geq N_{\max }$ $%
(r;f)$, and this later is larger than the more commonly used average with
respect to surface area measure \cite{rudin}, $\frac{1}{\sigma _{d}}%
\int_{S}N(r,\frac{1}{f_{\zeta }})d\sigma _{d}$, of the counting functions
for the slice functions of $f$. The main contribution in this note is the
observation that the employment of $N(r;f)$ leads to results completely
analogous with those in dimension $1$. Furthermore, these result are
dimension-free and sharp. It should be pointed out that, since $\log |f(z)|$
is subharmonic in $%
%TCIMACRO{\U{211d} }%
%BeginExpansion
\mathbb{R}
%EndExpansion
^{2d}$, then Hayman's theorem \cite{hay} on the maximum "modulus" of
subharmonic functions in higher dimensions, may be employed to obtain a
result about the size of the Riesz-mass of $f$, measured properly, in
comparison with the maximum modulus. But Hayman's result involves the
dimension $d.$

The order $\rho $ of the entire function $f$ of several complex variables is
defined as in the case when the dimension $d=1$, namely by $(1)$ above, and
it is immediate that if $f$ is of order $\rho $, then each slice function $%
f_{\zeta }$ is of order at most $\rho .$ Also, by Jensen's formula, 
\begin{equation*}
N(r,\frac{1}{f_{\zeta }})=\frac{1}{2\pi }\int_{-\pi }^{\pi }\log |f_{\zeta
}(re^{i\theta })|d\theta ,
\end{equation*}%
so that $N_{\max }(r)\leq \log M(r;f)$. In particular, the order of $N_{\max
}(r)$ is less than or equal the order of $f.$Jensen's formula also gives us
that $n_{\max }(t;f))\log 2\leq \log M(2r;f)$, so that $N(r;f)$ is a
continuous increasing function of order less than or equal the order of $f$.

\begin{theorem}
If $f$ $:%
%TCIMACRO{\U{2102} }%
%BeginExpansion
\mathbb{C}
%EndExpansion
^{d}\rightarrow 
%TCIMACRO{\U{2102} }%
%BeginExpansion
\mathbb{C}
%EndExpansion
$ is entire, $f(0)=1$, and its order $\rho \in (0,1)$, then 
\begin{equation}
\lim \sup_{r\rightarrow \infty }\frac{N(r;f)}{\log M(r;f)}\geq \lim
\sup_{r\rightarrow \infty }\frac{N_{\max }(r)}{\log M(r;f)}\geq \frac{\sin
\pi \rho }{\pi \rho },
\end{equation}%
and%
\begin{equation*}
\lim \sup_{r\rightarrow \infty }\frac{N(r,f)}{(m_{p})_{\max }(r;f)}\geq 
\frac{1}{\left\{ \frac{1}{2\pi }\int_{-\pi }^{\pi }|\psi (\theta
)|^{p}d\theta \right\} ^{1/p}},1\leq p<\infty .
\end{equation*}%
Furthermore, these inequalities are best possible.
\end{theorem}

\begin{proof}
If $X=\limsup_{r\rightarrow \infty }\frac{N_{\max }(r)}{\log M(r;f)}$, then
for $\epsilon >0$, there is an $R>0$ such that $N_{\max }(r)<(X+\epsilon
)\log M(r;f)$ for all $r\geq R$. Also, since $\log M(r;f)$ is of order $\rho 
$, there exists \cite{lev} a slowly varying function $L(r)$ and a sequence $%
r_{n}$ increasing to infinity such that 
\begin{eqnarray*}
(i)\text{ }\log M(r;f) &\leq &r^{\lambda }L(r),\text{ }(r>0); \\
(ii)\text{ }\log M(r_{n};f) &=&r_{n}^{\lambda }L(r_{n}).
\end{eqnarray*}%
Using $(i)$, and a well-known inequality for the maximum modulus of entire
functions of order less than $1$, we obtain, for sufficiently large n, that 
\begin{equation*}
\log M(r_{n};f)=\sup_{\zeta \in S}\log M(r_{n};f_{\zeta })\leq \sup_{\zeta
\in S}\int_{0}^{\infty }\frac{r_{n}N(t,\frac{1}{f_{\zeta }})}{(t+r_{n})^{2}}%
dt\leq \int_{0}^{\infty }\frac{r_{n}N_{\max }(t)}{(t+r_{n})^{2}}dt
\end{equation*}%
\begin{equation*}
\leq \int_{0}^{R}\frac{r_{n}N_{\max }(t)}{(t+r_{n})^{2}}dt+(X+\epsilon
)\int_{R}^{\infty }\frac{r_{n}\log M(t;f)}{(t+r_{n})^{2}}dt
\end{equation*}%
\begin{equation*}
\leq \frac{R}{R+r_{n}}\log M(R;f)+(X+\epsilon )\int_{R}^{\infty }\frac{%
r_{n}t^{\rho }L(t)}{(t+r_{n})^{2}}dt.
\end{equation*}%
Dividing both sides by $\log M(r_{n};f)$ and taking limits as $%
r_{n}\rightarrow \infty $, we obtain, using properties of regularly varying
functions, 
\begin{equation*}
1\leq (X+\epsilon )\lim \sup_{r_{n}\rightarrow \infty }\frac{1}{\log
M(r_{n};f)}\int_{R}^{\infty }\frac{r_{n}t^{\rho }L(t)}{(t+r_{n})^{2}}%
dt=\lim_{r\rightarrow \infty }\frac{1}{r^{\rho }L(r)}\int_{R}^{\infty }\frac{%
rt^{\rho }L(t)}{(t+r)^{2}}dt=\frac{\pi \rho }{\sin \pi \rho }.
\end{equation*}%
We conclude that 
\begin{equation*}
\lim \sup_{r\rightarrow \infty }\frac{N_{\max }(r)}{\log M(r;f)}\geq \frac{%
\sin \pi \rho }{\pi \rho }.
\end{equation*}
\end{proof}

\section{\protect\bigskip The Best Possible Character of the result}

It remains to demonstrate the sharpness of this result. For this purpose
take any canonical product of order $\rho \in (0,1)$, whose zeros are
regularly distributed along one ray. For example, we may take 
\begin{equation*}
P(z)=\prod_{n=1}^{\infty }(1+\frac{z}{r_{n}}),\text{ }(r_{n}=n^{1/\rho
},0<\rho <1).
\end{equation*}%
Now, for $z,\eta \in 
%TCIMACRO{\U{2102} }%
%BeginExpansion
\mathbb{C}
%EndExpansion
^{d}$, write $z=(z_{1},z_{2},...,z_{d})$, $\eta =(\eta _{1},\eta
_{2},...,\eta _{d})$ and recall that the inner product and the norm are
defined by $z\cdot \eta =\sum_{n=1}^{d}z_{j}\bar{\eta}_{j}$, and $\parallel
z\parallel =\sqrt{z\cdot z}$. Now let $g(z)=P(z\cdot \eta )$ where $\eta \in
S$ is fixed throughout. Thus 
\begin{equation*}
g(z)=\prod_{n=1}^{\infty }(1+\frac{z\cdot \eta }{r_{n}}),\text{ }%
(r_{n}>0,\eta \in S,z\in 
%TCIMACRO{\U{2102} }%
%BeginExpansion
\mathbb{C}
%EndExpansion
^{d},0<\rho <1).
\end{equation*}%
We shall calculate $M(r;g)$ and $N_{\max }(r).$ If $r>0,$ $z\in 
%TCIMACRO{\U{2102} }%
%BeginExpansion
\mathbb{C}
%EndExpansion
^{d}$, and $\parallel z\parallel =r,$then, since $|z\cdot \eta |\leq
\parallel z\parallel $, we have that 
\begin{equation*}
|g(z)|\leq \prod_{n=1}^{\infty }(1+\frac{\parallel z\parallel }{r_{n}})=P(r).
\end{equation*}%
On the other hand, $|g(r\eta )|=P(r)$. Hence $M(r;g)=P(r).$ Moving on to the
counting functions, if $\zeta \in S$, and $\zeta \cdot \eta =0$, then $N(r,%
\frac{1}{g_{\zeta }})$ is identically zero. If $\zeta \cdot \eta \neq 0$,
then the slice function may be written in the form $g_{\zeta }(z)=$ $%
\prod_{n=1}^{\infty }(1+\frac{z}{r_{n}/\zeta \cdot \eta })$ .Thus, if $z\in 
%TCIMACRO{\U{2102} }%
%BeginExpansion
\mathbb{C}
%EndExpansion
$ and $|z|=r$, then $N(r,\frac{1}{g_{\zeta }})$ is given by 
\begin{equation*}
N(r,\frac{1}{g_{\zeta }})=\sum_{r_{n}\leq r|\zeta \cdot \eta |}\log \frac{r}{%
r_{n}/|\zeta \cdot \eta |}\leq \sum_{r_{n}\leq r|\zeta \cdot \eta |}\log 
\frac{r}{r_{n}}\leq N(r,\frac{1}{g_{\eta }})=N(r,\frac{1}{P}),
\end{equation*}%
and it follows that $N_{\max }(r)=N(r,\frac{1}{P}).$

We now invoke a well-known abelian result to conclude that 
\begin{equation*}
\lim_{r\rightarrow \infty }\frac{N_{\max }(r)}{\log M(r;g)}%
=\lim_{r\rightarrow \infty }\frac{N(r,\frac{1}{P})}{\log P(r)}=\frac{\pi
\rho }{\sin \pi \rho }.
\end{equation*}%
This establishes the best possible character of the inequality in (2).

\subsection{The L$_{p}$ result}

In order to prove the second part of Theorem 1, we find it convenient to
introduce the function $v$ defined by 
\begin{equation*}
v(z)=\func{Re}\int_{0}^{\infty }\frac{z}{t(t+z)}n_{\max }(t)dt,z=re^{i\theta
},r\geq 0,|\theta |<\pi .
\end{equation*}%
Since the order of $n_{\max }$ is less than one, $v$ is well-defined and
harmonic in the plane slit along the negative real axis. If now we put for $%
r>0,0\leq \theta <\pi $%
\begin{equation*}
u(re^{i\theta })=\frac{1}{\pi }\int_{0}^{\theta }v(re^{i\omega })d\omega ,
\end{equation*}%
then%
\begin{equation*}
u(re^{i\theta })=\frac{1}{\pi }\int_{0}^{\infty }\frac{r\sin \theta }{%
t^{2}+2tr\cos \theta +r^{2}}N(t;f)dt,
\end{equation*}%
and it follows that $u$ vanishes on the positive real axis and has the
boundary value $N(r;f)$ on the negative real axis. Accordingly, we set $%
u(-r)=N(r;f),r\geq 0.$

Fix $\theta \in (0,\pi )$ and let $E$ be a set of Lebesgue measure $2\theta
. $Then for $\zeta \in S$, a well known computation [] gives us that 
\begin{equation*}
\frac{1}{\pi }\int_{E}\log |f_{\zeta }(re^{i\omega })|d\omega \leq \frac{1}{%
\pi }\int_{0}^{\infty }\frac{r\sin \theta }{t^{2}+2tr\cos \theta +r^{2}}N(t;%
\frac{1}{f_{\zeta }})dt
\end{equation*}%
\begin{equation*}
\leq \frac{1}{\pi }\int_{0}^{\infty }\frac{r\sin \theta }{t^{2}+2tr\cos
\theta +r^{2}}N(t;f)dt=\frac{1}{\pi }\int_{0}^{\theta }v(re^{i\omega
})d\omega .
\end{equation*}%
This inequality is also valid for $\theta =0$ and $\theta =\pi .$But $%
v(re^{i\theta })$ is an even function of $\theta $ and it is non-increasing
on $(0,\pi )$ as can be seen e.g. by writing it as a Stieltjes integral, and
so the last inequality implies that 
\begin{equation*}
(\log |f_{\zeta }|)^{\#}(re^{i\theta })=\sup_{|E|=2\theta }\frac{1}{\pi }%
\int_{E}\log |f_{\zeta }(re^{i\omega })|d\omega \leq \sup_{|E|=2\theta }%
\frac{1}{\pi }\int_{E}v(re^{i\omega })d\omega =v^{\#}(re^{i\theta }),
\end{equation*}%
where the sharp denotes the function introduced by Baernstein in his
definition of the star function.. Now using a result of Baernstein, we
conclude that 
\begin{equation*}
m_{p}(r,\log |f_{\zeta }|)\leq m_{p}(r,v),(1\leq p<\infty ),
\end{equation*}%
and hence that 
\begin{equation*}
\sup_{\zeta \in S}m_{p}(r,\log |f_{\zeta }|)\leq m_{p}(r,v),(1\leq p<\infty
).
\end{equation*}%
It remains to study the growth of $m_{p}(r,v)$ with respect to the function $%
N(r;f)$, and this will be carried out at a special sequence $r_{n}$
increasing to $\infty $ and satisfying%
\begin{equation*}
N(r_{n};f)\leq r^{\rho }L(r),r>0;N(r_{n};f)=r_{n}^{\rho }L(r_{n}).
\end{equation*}%
Here $L$ is a slowly varying function in the sense of Karamata i.e. $L$ is
positive and $L(\sigma r)/L(r)\rightarrow 1$, as $r\rightarrow \infty ,$ for
each $\sigma >1.$The existence of this sequence follows []from the fact that 
$N(r;f)$ is continuous, non-decreasing, and of order $\rho .$

\section{Extension To Vector-Valued Functions}

The following extension of the result in section 2 to vector-valued
functions is straightforward.

\begin{theorem}
Let $f$ $:%
%TCIMACRO{\U{2102} }%
%BeginExpansion
\mathbb{C}
%EndExpansion
^{m}\rightarrow 
%TCIMACRO{\U{2102} }%
%BeginExpansion
\mathbb{C}
%EndExpansion
^{n}$be an entire function of order $\rho \in (0,1)$ and satisfying $%
f(0)=(1,1,...,1)$. Write $f=(f^{1},f^{2},...,f^{n})$ and put 
\begin{equation*}
N_{\max }(r)=\max_{1\leq j\leq n}\{\sup_{\zeta \in S}N(r,\frac{1}{f_{\zeta
}^{j}})\},M_{\max }(r)=\max_{1\leq j\leq n}\{\sup_{\zeta \in S}M(r,f_{\zeta
}^{j})\},
\end{equation*}%
then 
\begin{equation}
\lim \sup_{r\rightarrow \infty }\frac{N_{\max }(r)}{\log M_{\max }(r;f)}\geq 
\frac{\sin \pi \rho }{\pi \rho }.
\end{equation}%
Furthermore, this inequality is best possible.
\end{theorem}

\bigskip

\end{document}